\documentclass[11pt]{article}
\usepackage{amssymb,amsmath,amsthm,amsfonts}
\usepackage[dvips]{graphicx, color}
\usepackage{subfigure}
\usepackage{latexsym}

\theoremstyle{plain}
\newtheorem{theorem}{Theorem}[section]
\newtheorem{lemma}{Lemma}[section]
\newtheorem{proposition}{Proposition}[section]

\numberwithin{equation}{section}
\theoremstyle{definition}
\newtheorem{definition}{Definition}

\newtheorem{remark}{Remark}[section]

\begin{document}

\title{\textbf{Local Null Controllability of a Chemotaxis System of
Parabolic-Elliptic Type}\thanks{%
This work is supported by the National Natural Science Foundation of China,
the National Basic Research Program of China (2011CB808002), the National
Research Foundation of South Africa, the National Science Foundation of
China (11201358), and the Fundamental Research Funds for the Central
Universities. }}
\date{April 18, 2013}
\author{Bao-Zhu Guo$^{a,c}$ ~and Liang Zhang$^{b,c}$\thanks{%
The corresponding author. Email: changleang@yahoo.com.cn} \\
%EndAName
$^a$\textit{Academy of Mathematics and Systems Science, Academia Sinica}\\
\textit{Beijing 100190, China} \\
$^b$\textit{Department of Mathematics,}\\
\textit{Wuhan University of Technology, Wuhan 430070, China}\\
$^c$\textit{School of Computational and Applied Mathematics}\\
\textit{University of the Witwatersrand, Wits 2050, Johannesburg, South
Africa}}
\maketitle

\begin{abstract}
In this paper, we are concerned with the controllability of a chemotaxis
system of parabolic-elliptic type. By linearizing the nonlinear system into
two separated linear equations to bypass the obstacle caused by the
nonlinear drift term, we establish the local null controllability of the
original nonlinear system. The approach is different from the usual way of
treating the coupled parabolic systems. \vspace{0.3cm}


\vspace{0.3cm}

\end{abstract}

\section{Introduction and main result}

In this paper, we are concerned with a controlled initial-boundary value
system of parabolic-elliptic type
\begin{equation}
\begin{cases}
\partial _{t}u=\nabla \cdot \left( \nabla u-\chi u\nabla v\right) +{\mathbf{1%
}}_{\omega }f & \mathrm{in}\ \Omega \times (0,T), \\
0=\Delta v-\gamma v+\delta u & \mathrm{in}\ \Omega \times (0,T), \\
{\partial }_{\nu }u=0,\partial _{\nu }v=0 & \mathrm{on}\ \partial \Omega
\times (0,T), \\
u(x,0)=u_{0}(x) & x\in \Omega,%
\end{cases}
\label{e}
\end{equation}
where and henceforth $u$ and $v$ are shorthands for states $u(x,t)$
and $v(x,t)$ at time spacial position $x\in \Omega$ and time $t\ge
0$, $\partial _{t}=\partial /\partial t$, $\partial _{\nu }=\partial
/\partial \nu $ is the directional derivative along the outward unit
normal vector $\nu $ on $\partial \Omega $, $\mathbf{1}_{\omega }$
stands for the characteristic function of $\omega $, $f$ is the
control function, $u_0$ is the initial value, $\chi ,\gamma $ and
$\delta $ are given positive constants.

\bigskip

The system (\ref{e}) without the control (i.e., $f\equiv 0$) is a simple
chemotaxis system, which was addressed by Keller and Segel \cite{keller} as
a model to describe the aggregation process in slime mold morphogenosis,
assuming that the cells emit directly the chemoattractant which is
immediately diffused. The unknown function $u$ then denotes the cell
density, whereas $v$ represents the concentration of the chemoattractant.
The validity of Keller-Segel's chemotaxis system is supported by some
experiments on the \textit{Escherichia Coli} bacteria and other interesting
physical interpretations. In fact, it has been extensively involved in many
medical and biological applications as well as some relevant areas such as
ecology and environment sciences. Moreover, because the model has a rich
structure from mathematical point of view, it deserves to be challenged.
Actually, some special but interesting cases have been studied such as the
aggregation, the blow-up of solutions, and the chemotactic collapse. Some
significant results have been achieved from different perspectives. We refer
to \cite{hillen} (also \cite{horstmann1}) for a survey where a quite
complete bibliography on the topic is included.

\bigskip\ In this paper, we study the Keller-Segel system from the
controllability point of view. We say that the system \eqref{e} is \textit{%
locally null controllable at time }$T$, if there exists a neighborhood of
the origin such that for any initial data $u_{0}$ belonging to this
neighborhood, the solution $(u,v)$ of \eqref{e} produced by corresponding
control function $f$ satisfies $u(x,T)=0$ for almost all $x\in \Omega$,
where the neighborhood and the control function space will be specified
later. Here, we consider the local null controllability instead of the exact
null controllability. The reason being interested is that the solutions of
Keller-Segel system may blow up in either finite time or infinite time,
which is shown in \cite{herrero} for the 3-d case. The 2-d case is even more
interesting and attractive. Actually, it has been found for the 2-d case
that the solution exists globally in finite time when the mass of the
initial value is less then a threshold value, while the solution will blow
up either in finite or in infinite time when the mass of initial value is
larger than the threshold value (see, e.g., \cite{jager}).

\bigskip

The study of the controllability for parabolic equations has been
thriving in the past  decade, see for instance \cite{barbu,fursikov}
and the references therein. Among them, a special interest is on the
controllability of coupled parabolic systems. For coupled systems,
the most practical situation is to impose the control force on one
equation, which has attracted intensive attention in the last few
years. We refer to a survey paper \cite{ammar0} wherein abundant
references are provided.

\bigskip

However, to the best of our knowledge, very few results are available to the
control problems of the Keller-Segel system \eqref{e} where a parabolic
equation is coupled with an elliptic equation through a drift term. Very
recently, the null controllability of a kind of nonlinear parabolic-elliptic
system of the following is considered in \cite{fernandez}:%
\begin{equation*}
\begin{cases}
\partial _{t}y-\Delta y=F(y,z)+{\mathbf{1}}_{\omega }f & \mathrm{in}\ \Omega
\times (0,T), \\
-\Delta z=f(y,z) & \mathrm{in}\ \Omega \times (0,T),%
\end{cases}%
\end{equation*}%
where $F(y,z)$ and $f(y,z)$ are nonlinear terms. In \cite{yagi1}, an
optimal control problem of the system \eqref{e} with the control to
be imposed on the second equation is considered. Our previous work
\cite{guo} is the first work that considers the local exact
controllability of a type of Keller-Segel system where a parabolic
equation couples with another parabolic equation. In the present
work, we attempt the controllability of system \eqref{e}, which is
probably the first work for this system. In the system \eqref{e},
since the drift term $-\chi \nabla \cdot (u\nabla v)$ destroys some
good properties of the diffusion operator which ensure the
regularity of the system, much more mathematical difficulties than
the aforementioned coupled parabolic systems are caused. These
include the regularity, the estimation of the \textquotedblleft
observability inequality\textquotedblright, and many others. The
usual way to establish the controllability of a nonlinear system is
to linearize the nonlinear system into some coupled linear ones.
Then, by the combination of the controllability result of the
linearized system and some fixed point results, one is able to
establish the controllability of the nonlinear system. We refer to
the typical works \cite{fernandez,guo} for the approach of this
kind. In this paper, we investigate, however, the controllability of
system \eqref{e} in a different way motivated intuitively by the
special mathematical structure of system \eqref{e}. We may decompose
this nonlinear system into two separated linear systems: one is a
controlled parabolic system, and another one is an irrelevant
elliptic equation. In such a way, we can bypass the obstacle caused
by the nonlinear drift term. This technique would be useful for
other coupled systems like drift-diffusion equations from the
semiconductor device.

\bigskip

Throughout the paper, $\Omega \subset {\mathbb{R}}^{N}(N\geq 1)$ stands for
a bounded domain with smooth boundary $\partial \Omega $, $\omega $ is a
nonempty open subset of $\Omega $, and $T>0$. $Q=\Omega \times (0,T)$, $%
\Sigma =\partial \Omega \times (0,T)$, and $Q_{\omega }=\omega \times (0,T)$%
. The norms of the usual Lebesgue function spaces $L^{p}(\Omega )$ and $%
L^{p}(Q)$ are denoted by $|\cdot |_{p}$ and $\Vert \cdot \Vert _{p},$
respectively. $W^{s,q}(\Omega )$, $W_{q}^{2,1}(Q)$, $C^{2,1}(\bar{Q })$ and $%
C^\alpha(\bar\Omega)$ ($s,\alpha \geq 0$, $1\leq q\leq \infty $) represent
the usual Sobolev spaces (see, e.g., \cite{lady}). When $q=2$, $H^{m}(\Omega
)=W^{m,2}(\Omega )$, $m\in {\mathbb{N}}$. In addition,
\begin{equation*}
W(0,T)=\left\{ y|y\in L^{2}(0,T;H^{1}(\Omega )),\partial _{t}y\in
L^{2}(0,T;H^{1}(\Omega )^{\ast })\right\}
\end{equation*}%
is equipped with the graph norm $\|y\|_{W(0,T)}=\|y\|_{L^2(0,T;H^1(%
\Omega))}+\|\partial_ty\|_{L^2(0,T;H^1(\Omega)^*)}$, where $H^{1}(\Omega
)^{\ast }$ denotes the dual space of $H^{1}(\Omega )$ and their duality
product is denoted by $\langle \cdot ,\cdot \rangle $. We also use $C$ to
denote a positive constant independent of $T$, which have different values
in different contexts.

\begin{definition}
\label{def} A pair of functions $(u,v),$ with $u\in W(0,T)\cap L^{\infty
}(Q) $ and $v\in L^{2}(0,T;H^{1}(\Omega ))\cap L^{\infty }(Q),$ is said to
be a \textit{weak solution of \eqref{e}} if for every $\varphi \in
L^{2}(0,T;H^{1}(\Omega ))$, the following identities hold true:
\begin{equation*}
\begin{array}{l}
\displaystyle \int_{0}^{T}\left\langle \partial _{t}u,\varphi \right\rangle
dt+\iint_{Q} \left[ (\nabla u-\chi u\nabla v)\cdot \nabla \varphi -{\mathbf{1%
}}_{\omega }f\varphi \right] dxdt =0, \cr\noalign{\vskip2mm}\displaystyle %
\iint_{Q}\left[ \nabla v\cdot \nabla \varphi +(\gamma v-\delta u)\varphi %
\right] dxdt =0.%
\end{array}%
\end{equation*}
\end{definition}

Now, we are in a position to state the main result of this paper.

\begin{theorem}
\label{th-non} Let $T>0$. For any initial value $u_0$ satisfying
\begin{equation}
\left\vert u_{0}\right\vert _{\infty }\leq e^{-c_{1}\left( 1+T+{\frac{1}{T}}%
\right) },  \label{intialdata}
\end{equation}%
where $c_1>0$ is a constant independent of $T$, there exists a control $f\in
L^{\infty }(Q_\omega)$ such that system \eqref{e} admits a solution $(u,v)$
satisfying
\begin{eqnarray*}
u\in C([0,T];L^2(\Omega))\cap L^{\infty }(Q),  v\in L^{2}(0,T;H^{1}(\Omega
))\cap L^{\infty }(Q),
\end{eqnarray*}
and $u(x,T)=0$ for almost all $x\in \Omega $.
\end{theorem}

\begin{remark}
Theorem \ref{th-non} gives an explicit representation of the initial data
with respect to time $T$. It particularly shows that the shorter of the
terminal time $T$ is, the smaller for the initial value for the null
controllability of system \eqref{e}. In addition, under the assumption of
Theorem \ref{th-non}, $\limsup_{t\rightarrow T^-}\|v(\cdot,t)\|_2=0$.
\end{remark}

We proceed as follows. In section 2, we give some preliminary results.
Section 3 is devoted to the null controllability of a scalar parabolic
equation, for which the $L^\infty$-control is obtained and its estimates
with respect to time $T$ are also established. The proof of Theorem \ref%
{th-non} is presented in section 4.

\section{Some results for linear equations}

In the sequel of the paper, we need some regularity results for linear
equations for both parabolic and elliptic types. We first consider the
well-posedness of the linear elliptic equation followed by
\begin{equation}
\begin{cases}
0=\Delta v-\gamma v+\delta \eta & \mathrm{in}\ \Omega , \\
{\partial }_{\nu }v=0 & \mathrm{on}\ \partial \Omega ,%
\end{cases}
\label{linear-e}
\end{equation}%
where $\gamma $ and $\delta $ are positive constants. The result of
Proposition \ref{lemmae} is brought from \cite{agmon}.

\begin{proposition}
\label{lemmae} For any $\eta \in L^{p}(\Omega ),p>1,$ Eq. %
\eqref{linear-e} admits a unique solution $v\in W^{2,p}(\Omega )$ with
\begin{equation*}
\left\Vert v\right\Vert _{W^{2,p}(\Omega )}\leq C\left\vert \eta \right\vert
_{p}.
\end{equation*}
\end{proposition}

Next, we consider the parabolic equation
\begin{equation}
\begin{cases}
\partial _{t}u=\Delta u-\nabla \cdot \left( Bu\right) +F & \mathrm{in}\ Q,
\\
{\partial }_{\nu }u=0 & \mathrm{on}\ \Sigma , \\
u(x,0)=u_{0}(x) & x\in \Omega .%
\end{cases}
\label{linear-p}
\end{equation}

\begin{proposition}
\label{lemmap} Let $B\in L^{\infty }(Q)^{N}$ with $B\cdot \nu =0$ on $\Sigma
$, $F\in L^{\infty }(Q)$, and $u_{0}\in L^{\infty }(\Omega )$. Then Eq.
\eqref{linear-p} admits a weak solution $u\in L^{\infty }(Q)$ with
\begin{equation}
\Vert u\Vert _{\infty }\leq e^{C\varrho _{0}}\left( |u_{0}|_{\infty }+\Vert
F\Vert _{\infty }\right) ,  \label{para0}
\end{equation}%
where $C=C(\Omega)$ is a positive constant depending on $\Omega$, and  $\varrho_0$ is given by
\begin{equation}
\varrho _{0}=(1+\left\Vert B\right\Vert _{\infty }^{2})(1+T).  \label{rho0}
\end{equation}
\end{proposition}

A similar inequality \eqref{para0} could be found in \cite{lady}, but here we
improve the estimate so that it depends on time $T$ explicitly. To this end,
we need the following lemma (see \cite[Lemma 5.6, p. 95]{lady}).

\begin{lemma}
\label{lady} Suppose that a sequence $Y_{s},s=0,1,2,\cdots $ of nonnegative
numbers satisfy a recursion relation%
\begin{equation*}
Y_{s+1}\leq cb^{s}Y_{s}^{1+\varepsilon },\text{ \ \ }s=0,1,2,\cdots
\end{equation*}%
with some positive constants $c,\varepsilon $ and $b\geq1.$ Then $%
Y_{s}\rightarrow 0$ as $s\rightarrow \infty$ provided that
\begin{equation*}
Y_{0}\leq c^{-\frac{1}{\varepsilon }}b^{-\frac{1}{\varepsilon ^{2}}}.
\end{equation*}
\end{lemma}

\noindent \textbf{Proof of Proposition \ref{lemmap}.}\quad Let $%
(u-k)_{+}=\max \{u-k,0\}$ and $A_{k}(t)=\mathrm{meas}\{x\in \Omega|u(x,t)>k\}$ for $t\in \lbrack 0,T]$, where%
\begin{equation}
k\geq K_{0}=\left\Vert F\right\Vert _{\infty }+\left\vert u_{0}\right\vert
_{\infty }.  \label{kkk0}
\end{equation}%
Multiplying by $(u-k)_{+}$ the both sides of \eqref{linear-p}, we get, by
integration by parts and the H\"{o}lder inequality, that%
\begin{equation*}
\begin{array}{ll}
\displaystyle\frac{d}{dt}\int_{\Omega }\left\vert (u-k)_{+}\right\vert
^{2}dx+\int_{\Omega }\left\vert \nabla (u-k)_{+}\right\vert ^{2}dx\cr%
\noalign{\vskip2mm}\displaystyle\leq \left\Vert B\right\Vert _{\infty
}^{2}\int_{A_{k}(t)}u^{2}dx+\int_{\Omega }\left\vert (u-k)_{+}\right\vert
^{2}dx+\int_{A_{k}(t)}F^{2}dx.\cr\noalign{\vskip2mm}\displaystyle\leq
(2\left\Vert B\right\Vert _{\infty }^{2}+1)\left( \int_{\Omega }\left\vert
(u-k)_{+}\right\vert ^{2}dx+\int_{A_{k}(t)}k^{2}dx\right) . &
\end{array}%
\end{equation*}%
From Gronwall's inequality, it follows that
\begin{equation}
\int_{\Omega }\left\vert (u-k)_{+}\right\vert
^{2}dx+\int_{0}^{t}\int_{\Omega }\left\vert \nabla (u-k)_{+}\right\vert
^{2}dxdt\leq e^{C\varrho _{0}}\int_{0}^{T}\int_{A_{k}(t)}k^{2}dxdt
\label{reg3}
\end{equation}%
for all $t\in \lbrack 0,T]$, where and in what follows $\varrho _{0}$ is
given by \eqref{rho0}. On the other hand, by Proposition I.3.2 of \cite%
{dibenedetto},
\begin{equation}
\left\Vert v\right\Vert _{\frac{2(N+2)}{N}}\leq C(\Omega )(1+T)^{\frac{N}{%
2(N+2)}}\left\Vert v\right\Vert _{V_{2}(Q)}  \label{guo1}
\end{equation}%
for any $v\in V_{2}(Q)$. Here $V_{2}(Q)=L^{\infty }(0,T;L^{2}(\Omega ))\cap
L^{2}(0,T;H^{1}(\Omega ))$ is endowed with its graph norm. Then, \eqref{reg3}
together with (\ref{guo1}) gives
\begin{equation}
\left\Vert (u-k)_{+}\right\Vert _{\frac{2(N+2)}{N}}\leq e^{C\varrho
_{0}}\int_{0}^{T}\int_{A_{k}(t)}k^{2}dxdt.  \label{reg4}
\end{equation}%
Let $\varphi (k)=$meas$\{(x,t)\in Q|u(x,t)>k\}.$ Then for any $h>k,$ we get,
from \eqref{reg4}, that%
\begin{equation*}
(h-k)^{2}\varphi (h)^{\frac{N}{N+2}}\leq \left\Vert (u-k)_{+}\right\Vert _{%
\frac{2(N+2)}{N}}^{2}\leq e^{C\varrho _{0}}\varphi (k)k^{2},
\end{equation*}%
which then gives
\begin{equation}
\varphi (h)\leq e^{C\varrho _{0}}\left( \frac{k}{h-k}\right) ^{\frac{2(N+2)}{%
N}}\varphi (k)^{\frac{N+2}{N}}.  \label{reg5}
\end{equation}%
Next, set $Y_{s}=\varphi (k_{s}),k_{s}=M(2-\frac{1}{2^{s}})$, and put $%
h=k_{s+1}$ and $k=k_{s}$ in \eqref{reg5} to get%
\begin{equation*}
Y_{s+1}\leq \tilde{c}4^{\tau }\left( 2^{\tau }\right)
^{s}Y_{s}^{1+\varepsilon },
\end{equation*}%
where $\tau =2(N+2)/N,$ $\varepsilon =2/N,$ and $\tilde{c}=e^{C\varrho _{0}}.
$ By Lemma \ref{lady}, we get $\varphi (2M)=0$ provided that
\begin{equation}
Y_{0}=\varphi (k_{0})=\varphi (M)\leq \left( \tilde{c}4^{\tau }\right) ^{-%
\frac{1}{\varepsilon }}(2^{\tau })^{-\frac{1}{e^{2}}},  \label{reg6}
\end{equation}%
for some positive real number $M$. To determine the value of $M$, let $m>1$
be an integer and $M=mK_{0}$, where $K_{0}\ $is given by \eqref{kkk0}. Put $%
h=M=mK_{0}$ and $k=K_{0}$ in \eqref{reg5} to get
\begin{equation}
\varphi (M)\leq \tilde{c}\left( \frac{1}{m-1}\right) ^{\tau }\varphi
(K_{0})^{1+\varepsilon }\leq \tilde{c}\left( \frac{1}{m-1}\right) ^{\tau
}T^{1+\varepsilon }\left( \text{meas}\Omega\right) ^{1+\varepsilon }.
\label{reg7}
\end{equation}%
Now, to get $\varphi (2M)=0$, we need to choose a proper $m$ such that %
\eqref{reg6} holds. Combining \eqref{reg6} and \eqref{reg7}, we only need
the integer $m$ to be such that
\begin{equation*}
m\geq 1+\tilde{c}^{\frac{1+\varepsilon }{\varepsilon \tau }}T^{\frac{%
1+\varepsilon }{\tau }}\left( \text{meas}\Omega\right) ^{\frac{1+\varepsilon }{%
\tau }}2^{\frac{2}{\varepsilon }+\frac{1}{\varepsilon ^{2}}}.
\end{equation*}%
Hence $\varphi (2M)=\varphi (2mK_{0})=0$ gives
\begin{equation*}
u\leq 2mK_{0}\leq e^{C\varrho _{0}}\left( \left\Vert F\right\Vert _{\infty
}+\left\vert u_{0}\right\vert _{\infty }\right) .
\end{equation*}%
In a similar argument, we can also get the other half part of \eqref{para0}
for $-u$. This completes the proof. \hfill $\Box $

\section{ Null controllability of a linear parabolic equation}

In this section, we consider the null controllability of the linear
parabolic equation
\begin{equation}
\begin{cases}
\partial _{t}u=\Delta u-\nabla \cdot \left( Bu\right) +\mathbf{1}_{\omega }f
& \mathrm{in}\ Q, \\
{\partial }_{\nu }u=0 & \mathrm{on}\ \Sigma , \\
u(x,0)=u_{0}(x) & x\in \Omega.%
\end{cases}
\label{linear-par0}
\end{equation}

\begin{theorem}
\label{th-lin} Let $T>0$, and $B\in L^{\infty }(Q)^N$ with $B\cdot \nu =0$
on $\Sigma $. For any $u_{0}\in L^{2}(\Omega )$, there exists a control $%
f\in L^{\infty }(Q_\omega)$ such that the solution $u$ of system %
\eqref{linear-par0} corresponding to $f$ satisfies $u\in W(0,T)$ and $%
u(x,T)=0$ for $x\in \Omega $ almost everywhere. Moreover, the control $f$
satisfies%
\begin{equation}
\left\Vert {\mathbf{1}}_\omega f\right\Vert _{\infty }\leq e^{C\kappa
}\left\vert u_{0}\right\vert _{2},  \label{fcontrol}
\end{equation}%
where
\begin{equation}
\kappa =(1+\Vert B\Vert _{\infty }^{2})(1+T)+{\frac{1}{T}}.  \label{kconst}
\end{equation}
\end{theorem}

To prove Theorem \ref{th-lin}, we need to establish a type of
``observability inequality'' for the following adjoint equation of %
\eqref{linear-par0}:
\begin{equation}
\begin{cases}
-\partial _{t}\phi =\Delta \phi +B\cdot\nabla \phi & \mathrm{in}\ Q, \\
{\partial }_{\nu }\phi =0 & \mathrm{on}\ \Sigma , \\
\phi (x,T)=\phi ^{T}(x) & x\in \Omega ,%
\end{cases}
\label{adjoint00}
\end{equation}%
where $\phi ^{T}\in L^{2}(\Omega ).$ By \cite[Lemma 1.1]{fursikov}, there is
a function $\beta \in C^{2}(\overline{\Omega })$ such that $\beta (x)>0$ for
all $x\in \Omega $ and $\beta |_{\partial \Omega }=0,\left\vert \nabla \beta
(x)\right\vert >0$ for all $x\in \overline{\Omega \setminus \omega }$. For $%
\lambda >0,$ set%
\begin{equation}
\varphi =\frac{e^{\lambda \beta }}{t(T-t)},\text{ \ }\alpha =\frac{%
e^{\lambda \beta }-e^{2\lambda \left\Vert \beta \right\Vert _{C(\overline{%
\Omega })}}}{t(T-t)}.  \label{alpha}
\end{equation}

We then have a Carleman inequality stated in Lemma \ref{carleman} (see \cite%
{fursikov}).

\begin{lemma}
\label{carleman}There exists a constant $\lambda _{0}=\lambda _{0}(\Omega
,\omega )>1$ such that for all $\lambda \geq \lambda _{0}$ and $s\geq \gamma
(\lambda )(T+T^{2})$,
\begin{eqnarray}
&&\iint_{Q}\left[ s\varphi |\nabla y|^{2}+(s\varphi )^{3}|y|^{2}\right]
e^{2s\alpha }\ dxdt  \notag \\
&\leq &C\iint_{Q}e^{2s\alpha }|\partial _{t}y\pm\Delta y|^{2}\
dxdt+\iint_{Q_{\omega }}(s\varphi )^{3}e^{2s\alpha }|y|^{2}\ dxdt
\label{carle0}
\end{eqnarray}%
for all $y\in X= \{\xi \in C^{2,1}(\bar{Q})|\partial _{\nu }\xi =0$ on $%
\Sigma \},$ where $\gamma (\lambda )$ is given by
\begin{equation}
\gamma (\lambda )=e^{2\lambda \left\Vert \beta \right\Vert _{C(\overline{%
\Omega })}}.  \label{gamma}
\end{equation}
\end{lemma}

Proposition \ref{pro-obser} is an ``observability inequality'' for the
adjoint equation \eqref{adjoint00}.

\begin{proposition}
\label{pro-obser} Let $\delta _{0}\in (1,2)$. Then there exist positive
constants $\lambda $ and $s$ such that for all $T>0,\phi ^{T}\in
L^{2}(\Omega ),$ the solution $\phi $ of system \eqref{adjoint00} satisfies%
\begin{equation}
\left\vert \phi (\cdot ,0)\right\vert _{2}^{2}\leq e^{C\kappa
}\iint_{Q_{\omega }}e^{\delta _{0}s\alpha }\left\vert \phi \right\vert
^{2}dxdt,  \label{observability}
\end{equation}%
where $\kappa $ is given by \eqref{kconst}.
\end{proposition}

\noindent \textit{Proof.}\quad First, by Lemma \ref{carleman}, there exists
a positive constant $\lambda _{1}=C_1(\Omega ,\omega )(1+\left\Vert
B\right\Vert _{\infty }^{2})$ satisfying $\gamma (\lambda _{1})\geq \lambda
_{1}>1$ such that for any $\lambda \geq \lambda _{1},s\geq \gamma (\lambda
)(T+T^{2})$ and $\phi ^{T}\in L^{2}(\Omega ),$ the associated solution $\phi
$ to \eqref{adjoint00} satisfies%
\begin{equation}
\iint_{Q}(s\varphi )^{3}|\phi |^{2}e^{2s\alpha }\ dxdt\leq C\iint_{Q_{\omega
}}(s\varphi )^{3}\left\vert \phi \right\vert ^{2}e^{2s\alpha }dxdt,
\label{carle1}
\end{equation}%
where $C=C(\Omega ,\omega )$, and $\gamma (\lambda _{1})$ and $\gamma
(\lambda )$ are given by \eqref{gamma}. By \eqref{adjoint00},
\begin{equation*}
\frac{d}{dt}\left( e^{\left\Vert B\right\Vert _{\infty }^{2}t}\left\vert
\phi \right\vert _{2}^{2}\right) \geq 0.
\end{equation*}%
This gives, for any $t\in (0,T]$, that
\begin{equation}
\left\vert \phi (\cdot ,0)\right\vert _{2}^{2}\leq e^{\left\Vert
B\right\Vert _{\infty }^{2}T}\left\vert \phi (\cdot ,t)\right\vert _{2}^{2}
\label{aj1}
\end{equation}%
for all $t\in (0,T]$. Integrate both sides of \eqref{aj1} over $\left[
T/4,3T/4\right] $ to give
\begin{equation}
\left\vert \phi (\cdot ,0)\right\vert _{2}^{2}\leq \frac{2}{T}e^{\left\Vert
B\right\Vert _{\infty }^{2}T}\int_{\frac{T}{4}}^{\frac{3T}{4}}\int_{\Omega
}\left\vert \phi \right\vert ^{2}dxdt.  \label{aj2}
\end{equation}%
Since $(s\varphi )^{-3}e^{-2s\alpha }\leq e^{Cs/T^{2}}$ in $\Omega \times
\lbrack T/4,3T/4],$ inequality \eqref{observability} then follows from %
\eqref{carle1} and \eqref{aj2}, with $\lambda $ and $s$ taken as $\lambda
=C(1+\left\Vert B\right\Vert _{\infty }^{2})>\lambda _{1}$ and $%
s=C(1+\left\Vert B\right\Vert _{\infty }^{2})(T+T^{2}).$ \hfill $\Box $

\bigskip

\noindent \textbf{Proof of Theorem \ref{th-lin}.} Let $\varepsilon >0.$ We
set%
\begin{equation*}
J_{\varepsilon }(u,f)=\frac{1}{2}\iint_{Q_{\omega }}\left\vert f\right\vert
^{2}e^{-\delta _{0}s\alpha }dxdt+\frac{1}{2\varepsilon }\int_{\Omega
}\left\vert u(x,T)\right\vert ^{2}dx
\end{equation*}%
and consider the extremal problem
\begin{equation*}
\inf_{(u,f)\in \mathcal{U}}J_{\varepsilon }(u,f),
\end{equation*}%
where $\mathcal{U}$ is the totality of $(u,f)\in W(0,T)\times L^{2}(Q)$
solving \eqref{linear-par0}. The existence of an optimal pair $%
(f_{\varepsilon },u_{\varepsilon })$ to the above extremal problem follows
from a standard argument. By the maximum principle (see \cite{barbu1}), we
get the optimality system for this problem as follows:
\begin{eqnarray}
&&%
\begin{cases}
-\partial _{t}\phi _{\varepsilon }=\Delta \phi _{\varepsilon }+B\cdot\nabla
\phi _{\varepsilon } & \mathrm{in}\ Q, \\
{\partial }_{\nu }\phi _{\varepsilon }=0 & \mathrm{on}\ \Sigma , \\
\phi _{\varepsilon }(x,T)=-\frac{1}{\varepsilon }u_{\varepsilon }(x,T) &
x\in \Omega ;%
\end{cases}
\label{adjointe} \\
&&%
\begin{cases}
\partial _{t}u_{\varepsilon }=\Delta u_{\varepsilon }-\nabla \cdot \left(
Bu_{\varepsilon }\right) +\mathbf{1}_{\omega }f_{\varepsilon } & \mathrm{in}%
\ Q, \\
{\partial }_{\nu }u_{\varepsilon }=0 & \mathrm{on}\ \Sigma , \\
u_{\varepsilon }(x,0)=u_{0}(x) & x\in \Omega ;%
\end{cases}
\label{elineare} \\
&&f_{\varepsilon }-\mathbf{1}_{\omega }\phi _{\varepsilon }e^{\delta
_{0}s\alpha }=0.  \label{fcontrols}
\end{eqnarray}%
Note that we can choose $s$ and $\lambda $ such that the observability
inequality \eqref{observability} and
\begin{equation}
\omega (\lambda )=e^{-\lambda \left\Vert \beta \right\Vert _{C(\overline{%
\Omega })}}<\delta _{0}-1  \label{eta}
\end{equation}%
hold. Furthermore, by \eqref{adjointe}, \eqref{elineare}, \eqref{fcontrols},
and Proposition \ref{pro-obser}, we get%
\begin{equation}
\iint_{Q_{\omega }}\left\vert \phi _{\varepsilon }\right\vert ^{2}e^{\delta
_{0}s\alpha }dxdt+\frac{1}{\varepsilon }\int_{\Omega }\left\vert
u_\varepsilon(x,T)\right\vert ^{2}dx\leq e^{C\kappa }\left\vert u_{0}\right\vert
_{2}^{2}.  \label{null0}
\end{equation}%
Eq. \eqref{fcontrols} and \eqref{null0} then lead to $\left\Vert{%
\mathbf{1}}_\omega f_{\varepsilon }\right\Vert _{2}\leq e^{C\kappa
}\left\vert u_{0}\right\vert _{2},$ which means that the controls $%
f_\varepsilon$ can be taken in  $L^{2}$ space.

Now, we show that $f_{\varepsilon}$ can actually be taken in   $L^{\infty }$ space.
To this purpose, we apply a so-called bootstrap method in \cite{zhang}
(see also \cite{barbu}). Firstly, set $\alpha _{0}=\min_{\overline{\Omega }%
}\alpha$. Then the following inequalities can be easily verified:
\begin{equation}
\alpha _{0}\leq \alpha \leq \frac{\alpha _{0}}{1+\omega (\lambda )}<0,
\label{alpha0}
\end{equation}%
where $\omega (\lambda )$ is defined by \eqref{eta}. Secondly, let $\tau $
be a sufficiently small positive constant and let $\left\{ \tau _{j}\right\}
_{j=0}^{M}$ be a finite increasing sequence such that $0<\tau _{j}<\tau
,j=0,1,\ldots ,M,\tau _{M}=\tau.$ Let $\left\{ p_{i}\right\} _{i=0}^{M}$ be
another finite increasing sequence such that $p_{0}=2,p_{M}=\infty $, and
\begin{equation}
-\left( \frac{N}{2}+1\right) \left( \frac{1}{p_{i}}-\frac{1}{p_{i+1}}\right)
+1>\frac{1}{2},\ i=0,1,\ldots ,M-1.  \label{pipi}
\end{equation}%
For each $i,i=0,1,\ldots ,M,$ define%
\begin{eqnarray*}
z_{i}(x,t) &=&e^{(s+\tau _{i})\alpha _{0}}\phi _{\varepsilon }(x,T-t), \\
F_{i}(x,t) &=&\left[ \partial _{t}(e^{(s+\tau _{i})\alpha _{0}})\right] \phi
_{\varepsilon }(x,T-t), \\
\tilde{B}(x,t) &=&B(x,T-t).
\end{eqnarray*}%
Then, the adjoint equation \eqref{adjointe} is transformed into the
following initial-boundary problem for every $i=0,1,\ldots,M$:
\begin{equation}
\begin{cases}
\partial _{t}z_{i}-\Delta z_{i}=\tilde{B}\cdot\nabla z_{i}+F_{i} & \mathrm{in%
}\ Q, \\
{\partial }_{\nu }z_{i}=0 & \mathrm{on}\ \Sigma , \\
z_{i}(x,0)=0 & x\in \Omega .%
\end{cases}
\label{adjointee}
\end{equation}%
Let $\{S(t)\}_{t\geq 0}$ be the semigroup generated by the Laplace operator
with Neumann boundary condition. It follows that (see, e.g., \cite{arendt})%
\begin{equation}
\left\vert S(t)u\right\vert _{q}\leq Cm(t)^{-\frac{N}{2}\left( \frac{1}{p}-%
\frac{1}{q}\right) }\left\vert u\right\vert _{p}  \label{lp-lq}
\end{equation}%
for all $u\in L^{p}(\Omega ),$ $t>0$, and $1<p\leq q\leq \infty $, where $%
m(t)=\min \{1,t\}$. Note that the solution $z_{i}$ of \eqref{adjointee} can
be represented as%
\begin{equation}
z_{i}(\cdot ,t)=\int_{0}^{t}S(t-s)\left( \tilde{B}\cdot\nabla
z_{i}+F_{i}\right) (\cdot ,s)ds, \; i=0,1,\ldots ,M.  \label{rhoi}
\end{equation}%
Applying the estimates \eqref{lp-lq} to \eqref{rhoi}, we have
\begin{equation}  \label{guo3}
\left\vert z_{i}(\cdot ,t)\right\vert _{p_{i}}\leq C\int_{0}^{t}m(t-s)^{-%
\frac{N}{2}\left( \frac{1}{p_{i-1}}-\frac{1}{p_{i}}\right) }\left\vert
\left( \tilde{B}\cdot\nabla z_{i}+F_{i}\right) (\cdot ,s)\right\vert
_{p_{i-1}}ds,\; i=0,1,\ldots ,M.
\end{equation}%
\newline
With \eqref{pipi}, we apply Young's convolution inequality (see, e.g. \cite%
{arendt}) to the right-hand side of (\ref{guo3}) to get
\begin{equation}
\left\Vert z_{i}\right\Vert _{p_{i}}\leq e^{C(1+T)}\left( \left\Vert
B\right\Vert _{\infty }\left\Vert \nabla z_{i}\right\Vert
_{p_{i-1}}+\left\Vert F_{i}\right\Vert _{p_{i-1}}\right) .  \label{zzeta}
\end{equation}%
On the other hand, by a standard energy estimate applied to \eqref{adjointee},   we can get the
following $L^{p_{i-1}}$-estimate for $z_{i}$:
\begin{equation}
\left\Vert z_{i}\right\Vert _{p_{i-1}}+\left\Vert \nabla z_{i}\right\Vert
_{p_{i-1}}\leq e^{C\kappa }\left\Vert F_{i}\right\Vert _{p_{i-1}}.
\label{energg}
\end{equation}%
By the definition of $F_{i}$ and \eqref{alpha0}, we have
\begin{equation}
\left\Vert F_{i}\right\Vert _{p_{i-1}}\leq CT\left\Vert z_{i-1}\right\Vert
_{p_{i-1}}\text{ }.  \label{ghi}
\end{equation}%
Combining \eqref{zzeta}, \eqref{energg}, and \eqref{ghi} then, we obtain
\begin{equation*}
\left\Vert z_{i}\right\Vert _{p_{i}}\leq e^{C\kappa }\left\Vert
z_{i-1}\right\Vert _{p_{i-1}}.
\end{equation*}%
This iteration inequality from $0$ to $M$ produces
\begin{equation}
\left\Vert z_{M}\right\Vert _{p_{M}}\leq e^{C\kappa }\left\Vert
z_{0}\right\Vert _{2}.  \label{znn}
\end{equation}%
Since $p_{M}=\infty ,$ it follow from the definition of $z_{0}$, %
\eqref{null0}, and \eqref{znn} that $\left\Vert z_{M}\right\Vert _{\infty
}\leq e^{C\kappa }\left\Vert u_{0}\right\Vert _{2}; $ that is, $\left\Vert
\phi _{\varepsilon }e^{(s+\tau )\alpha _{0}}\right\Vert _{\infty }\leq
e^{C\kappa }\left\Vert u_{0}\right\Vert _{2}.$ By \eqref{fcontrols}, we get
\begin{equation}  \label{guo4}
\left\Vert e^{\left[ -s(\delta _{0}-1-\omega (\lambda ))+\tau (1+\omega
(\lambda ))\right] \alpha }{\mathbf{1}}_\omega f_{\varepsilon }\right\Vert
_{\infty }\leq e^{C\kappa }\left\Vert u_{0}\right\Vert _{2},
\end{equation}
where $\omega (\lambda )$ is given by \eqref{eta}. By choosing $\tau $ small
enough so that $-s\left( \delta _{0}-1-\omega (\lambda )\right) +\tau
(1+\omega (\lambda ))<0,$ we conclude from (\ref{guo4}) that
\begin{equation}
\left\Vert{\mathbf{1}}_\omega f_{\varepsilon }\right\Vert _{\infty }\leq
e^{C\kappa }\left\Vert u_{0}\right\Vert _{2}.  \label{fconte}
\end{equation}
This shows that the controls $f_\varepsilon$ can be taken in   $L^{\infty }$
space.

Finally, by \eqref{fconte}, we can extract a subsequences of $%
\{f_{\varepsilon }\}_{\varepsilon \geq 0}$, still denoted by itself, such
that ${\mathbf{1}}_{\omega }f_{\varepsilon }\rightarrow {\mathbf{1}}_{\omega
}f\in L^{\infty }(Q)$ weakly in $L^{2}(Q)$ as $\varepsilon \rightarrow 0$.
Denote by $u_{\varepsilon }$ the solution to the system \eqref{elineare} associated to $%
f_{\varepsilon }$. By virtue of Proposition \ref{lemmap}, $\{u_{\varepsilon
}\}_{\varepsilon \geq 0}$ is uniformly bounded in $W(0,T)$. Thus, we can
extract a subsequence of $\{u_{\varepsilon }\}_{\varepsilon \geq 0}$, still
denoted by itself, such that $u_{\varepsilon }\rightarrow u$ weakly in $%
W(0,T)$ for $u\in W(0,T)\subset C([0,T];L^{2}(\Omega ))$. Such a $u$ is the
weak solution of \eqref{linear-par0} corresponding to $f$. In addition, by \eqref{null0}, $u(x,T)=0$ almost everywhere in $\Omega $%
. This completes the proof. \hfill $\Box $

\section{Proof of Theorem \protect\ref{th-non}}

Let $K=\left\{ \xi \in L^{\infty }(Q)|\left\Vert \xi \right\Vert _{\infty
}\leq 1\right\} \cap L^{\infty }(0,T;L^{p}(\Omega ))\subset L^{2}(Q)$, $%
p>\max \{N,2\}$. For every $\xi \in K$, consider the following two linear
equations:
\begin{equation}
\begin{cases}
0=\Delta v(\cdot ,t)-\gamma v(\cdot ,t)+\delta \xi (\cdot ,t) & \mathrm{in}\
\Omega , \\
{\partial }_{\nu }v(\cdot ,t)=0 & \mathrm{on}\ \partial \Omega%
\end{cases}
\label{lin-e}
\end{equation}%
for almost every $t\in \lbrack 0,T]$ and%
\begin{equation}
\begin{cases}
\partial _{t}u=\Delta u-\nabla \cdot \left( Bu\right) +\mathbf{1}_{\omega }f
& \mathrm{in}\ Q, \\
{\partial }_{\nu }u=0 & \mathrm{on}\ \Sigma , \\
u(x,0)=u_{0}(x) & x\in \Omega ,%
\end{cases}
\label{lin-p}
\end{equation}%
where $B=\chi \nabla v^{\xi }$. In what follows, we denote by
shorthand $v^{\xi }=v^{\xi }(\cdot ,t)$ the unique solution of
equation \eqref{lin-e} corresponding to $\xi (\cdot ,t)$ for $t\in
\lbrack 0,T]$. First, by Proposition \ref{lemmae}, we see that
$v^{\xi }\in L^{\infty }(0,T;W^{2,p}(\Omega ))$ for  $p>\max\{N,2\}$
provided that $\xi \in K.$ Hence, the
embedding theory between Sobolev spaces for $p>N$ (see, e.g., \cite%
{dibenedetto}) then implies
\begin{equation*}
B=\chi \nabla v^{\xi }\in L^{\infty }(0,T;W^{1,p}(\Omega
))^{N}\subset
L^{\infty }(0,T;C(\bar{\Omega}))^{N}\text{ with }B\cdot \nu =0\text{ on }%
\Sigma ,
\end{equation*}%
and in addition,
\begin{equation}
\left\Vert B\right\Vert _{\infty }=\chi \left\Vert \nabla v^{\xi
}\right\Vert _{L^{\infty }(0,T;W^{1,p}(\Omega ))}\leq C\Vert \xi
\Vert _{\infty }\leq C.  \label{bbbb}
\end{equation}%
Thus, we can define a linear continuous operator $\Phi $ from $K\ $to $%
L^{\infty }(0,T;C(\bar{\Omega}))^{N}\subset L^{\infty }(Q)^{N}$ by
\begin{equation*}
\Phi (\xi )=B=\chi \nabla v^{\xi },\forall \;\xi \in K.
\end{equation*}%
Second, by Theorem \ref{th-lin}, we see that for each $B\in L^{\infty
}(Q)^{N}$ with $B\cdot \nu =0$ on $\Sigma $, there exists a pair $\left(
u,f\right) \in L^{2}(Q)\times L^{\infty }(Q)$ that solves system %
\eqref{lin-p} with $u(x,T)=0$ for almost all $x\in \Omega $. Moreover, the
control $f$ satisfies \eqref{fcontrol}. By \eqref{bbbb},
\begin{equation}
\left\Vert {\mathbf{1}}_{\omega }f\right\Vert _{\infty }\leq e^{C\kappa
_{0}}\left\vert u_{0}\right\vert _{2},  \label{fcontrolin0}
\end{equation}%
where
\begin{equation}
\kappa _{0}=  1+T+\frac{1}{T} .  \label{k0}
\end{equation}and in the sequel,  $C$ is a positive constant independent of time $T$.
By \eqref{kkk0} in Proposition \ref{lemmap} and \eqref{fcontrolin0}, we have
the following estimate:%
\begin{equation}
\left\Vert u\right\Vert _{W(0,T)}+\left\Vert u\right\Vert _{\infty }\leq
e^{C\kappa _{0}}\left\vert u_{0}\right\vert _{\infty }.  \label{rproof0}
\end{equation}%
We then define a multi-valued mapping $\Psi: L^{\infty }(Q)^{N} \to
2^{L^{2}(Q)}$ by
\begin{equation*}
\Psi (B)=\left\{
\begin{tabular}{ll}
$u\in L^{2}(Q)$ & $\left\vert
\begin{array}{c}
\exists f\in L^{\infty }(Q_{\omega })\text{ satisfying \eqref{fcontrolin0}
such that }u\text{ is the solution } \\
\text{ of \eqref{lin-p} corresponding to }f\text{ and }B,\text{ and }u(x,T)=0%
\text{ a.e. in }\Omega%
\end{array}%
\right. $%
\end{tabular}%
\right\}
\end{equation*}%
where $2^{L^{2}(Q)}$ stands for all subsets of $L^{2}(Q)$.
 Since both operators $\Phi $ and $\Psi $ are well defined,
which is
guaranteed by Proposition \ref{lemmae} and Theorem \ref{th-lin}, we let%
\begin{equation}
\Lambda =\Psi \circ \Phi :K\subset L^{2}(Q)\rightarrow 2^{L^{2}(Q)}.
\label{guo5}
\end{equation}%
Now, we apply Kakutani's fixed point theorem (see \cite[p.7]{barbu1}) to the
map $\Lambda $ to prove Theorem \ref{th-non}. Indeed, it is clear that $K$
is a convex subset of $L^{2}(Q).$ By Proposition \ref{lemmae} and Theorem %
\ref{th-lin} again, for any $\xi \in K,$ $\Lambda (\xi )$ is nonempty and it
is also convex due to the linearity of the equations. Moreover, from %
\eqref{rproof0}, it follows that for each $\xi \in K$, $\Lambda (\xi )$ is
bounded in $W(0,T)$, and hence a compact subset of $L^{2}(Q)$ according to
the Aubin-Lions lemma (see \cite[p.17]{barbu1}).

We claim that $\Lambda $ is upper semi-continuous. Indeed, let $\left\{ \xi
_{n}\right\} _{n=1}^{\infty }$ be a sequence of functions in $K$ such that
\begin{equation}
\xi _{n}\rightarrow \xi \text{ strongly in }L^{2}(Q)\text{ as }n\rightarrow
\infty .  \label{etannn}
\end{equation}%
For every $n$, let $B_{n}=\Phi (\xi _{n})=\chi \nabla v_{n}$ and take $%
u_{n}\in \Lambda (\xi _{n})=\Psi (B_{n})$, where $v_{n}$ solves
\begin{equation}
\begin{cases}
0=\Delta v_{n}(\cdot ,t)-\gamma v_{n}(\cdot ,t)+\delta \xi _{n}(\cdot ,t) &
\mathrm{in}\ \Omega , \\
{\partial }_{\nu }v_{n}(\cdot ,t)=0 & \mathrm{on}\ \partial \Omega%
\end{cases}
\label{lin-ee}
\end{equation}%
for almost every $t\in \lbrack 0,T]$ and $u_{n}$ solves%
\begin{equation}
\begin{cases}
\partial _{t}u_{n}=\Delta u_{n}-\nabla \cdot \left( B_{n}u_{n}\right) +%
\mathbf{1}_{\omega }f_{n} & \mathrm{in}\ Q, \\
{\partial }_{\nu }u_{n}=0 & \mathrm{on}\ \Sigma , \\
u_{n}(x,0)=u_{0}(x) & x\in \Omega ,%
\end{cases}
\label{lin-pp}
\end{equation}%
with $u_{n}(x,T)=0$ for almost all $x\in \Omega $. Moreover, the control $%
f_{n}$ satisfies%
\begin{equation}
\left\Vert {\mathbf{1}}_{\omega }f_{n}\right\Vert _{\infty }\leq e^{C\kappa
_{0}}\left\vert u_{0}\right\vert _{2}.  \label{fcn}
\end{equation}%
To show that $\Lambda $ is upper semi-continuous, it suffices to
prove that there exist a subsequence of $\{u_{n}\}_{n=1}^{\infty }$
such that it converges strongly to an element of $\Lambda (\eta )$
in $L^{2}(Q)$ topology.

In what follows,  we do not distinguish the sequence and its
subsequence by abuse of notation. First, the estimate \eqref{fcn}
enables us to obtain a function $f\in L^{\infty }(Q)$ and a
subsequence of $\{f_{n}\}_{n=1}^{\infty
},$ such that%
\begin{equation}
{\mathbf{1}}_{\omega }f_{n}\rightarrow {\mathbf{1}}_{\omega }f\text{ weakly
in }L^{2}(Q)\text{; weakly}^{\ast }\text{ in }L^{\infty }(Q)\hbox{ as
}n\rightarrow \infty .  \label{fcn1}
\end{equation}%
By \eqref{fcn} and Proposition \ref{lemmap}, $u_n$ satisfies \eqref{rproof0}; that is
\begin{equation}
\left\Vert u_{n}\right\Vert _{W(0,T)}+\left\Vert u_{n}\right\Vert _{\infty
}\leq e^{C\kappa _{0}}\left\vert u_{0}\right\vert _{2}.  \label{ynun}
\end{equation}%
Applying the Aubin-Lions lemma again, we get a $u\in W(0,T)\cap L^{\infty }(Q)$
and a subsequence of $\{u_{n}\}_{n=1}^{\infty }$ such that%
\begin{equation}
u_{n}\rightarrow u\text{ weakly in }W(0,T)\text{; strongly in }L^{2}(Q)\text{%
, as }n\rightarrow \infty .  \label{unli}
\end{equation}%
Furthermore, by the strong convergence of $\{u_{n}\}_{n=1}^\infty$ in $L^{2}(Q)$, we can
extract a subsequence of $\{u_{{n}}\}_{n=1}^{\infty }$ (see \cite[Lemma 2.1,
p. 72]{lady}) such that%
\begin{equation}
u_{n}\rightarrow u\text{ almost everywhere in }Q\text{ as }n\rightarrow
\infty .  \label{unli1}
\end{equation}%
On the other hand, by Proposition \ref{lemmae}, for each $n$ and $p>1$, it holds that
\begin{equation}
\left\Vert v_{n}(\cdot ,t)\right\Vert _{W^{2,p}(\Omega )}\leq C\left\Vert
\xi _{n}(\cdot ,t)\right\Vert _{L^{p}(\Omega )}\leq C, \;\hbox{for almost every}\; t\in
[0,T]   \label{ynvn}
\end{equation}%
 Thus,  $B_{n}$ satisfies %
\eqref{bbbb}; that is
\begin{equation}
\left\Vert B_{n}\right\Vert _{\infty }\leq C\left\Vert \nabla
v_{n}\right\Vert _{L^{\infty }(0,T;W^{1,p}(\Omega ))}\leq C\Vert \xi
_{n}\Vert _{\infty }\leq C.  \label{ynbn}
\end{equation}%
Furthermore, let $v$ be the unique solution of \eqref{lin-e} corresponding
to $\xi .$ Then, by the linearity of Eq. \eqref{lin-e} and by %
\eqref{ynvn},
\begin{equation*}
\left\Vert v_{n}(\cdot ,t)-v(\cdot ,t)\right\Vert _{H^{2}(\Omega )}\leq
C\left\vert \xi _{n}(\cdot ,t)-\xi (\cdot ,t)\right\vert _{2}, \; \hbox{for almost every}\;
t\in [0,T].
\end{equation*}%
Since $\xi _{n}\rightarrow \eta $ strongly in $L^{2}(Q)$ in condition \eqref{etannn}, it follows that
\begin{equation*}
\left\Vert v_{n}-v\right\Vert _{L^{2}(0,T;H^{2}(\Omega ))}\leq C\left\Vert
\xi _{n}-\xi \right\Vert _{2}\rightarrow 0,
\end{equation*}%
which implies that
\begin{equation}
B_{n}=\chi \nabla v_{n}\rightarrow \chi \nabla v=B\ \text{strongly in }%
L^{2}(Q).  \label{ynbn1}
\end{equation}%
By \eqref{unli} and \eqref{ynbn},  the sequence $%
\{B_{n}u_{n}\}_{n=1}^{\infty }$ is bounded in $L^{2}(Q)$, whence
there is a subsequence  such that $B_{n}u_{n}\rightarrow \pi $
weakly in $L^{2}(Q)$ as $n\rightarrow \infty .$ By \eqref{ynbn1},
there is a subsequence of $\{B_{n}\}_{n=1}^\infty$ such that
$B_{n}\rightarrow B$ almost everywhere in $Q$. This together with
\eqref{unli1} implies that $B_{n}u_{n}\rightarrow Bu$ almost
everywhere in $Q $. Therefore $\pi =Bu$,  and
\begin{equation}
B_{n}u_{n}\rightarrow Bu\text{ weakly in }L^{2}(Q)\text{ as }n\rightarrow
\infty .  \label{bnun}
\end{equation}%
Now, by \eqref{fcn1}, \eqref{unli}, and \eqref{bnun}, we can pass to the
limit as $n\rightarrow \infty $ in \eqref{lin-pp} to get that $u$ is a weak
solution of \eqref{lin-p} corresponding to $B$ and $f$ in the sense of
Definition \ref{def}. Now, it only need to show that $u\in \Lambda (\xi )=\Psi (B)$.
Actually, let $U_{n}=u_{n}-u$ and $F_{n}=\mathbf{1}_{\omega }(f_{n}-f)$.
Then $U_{n}$ solves the system%
\begin{equation}
\begin{cases}
\partial _{t}U_{n}=\Delta U_{n}-\nabla \cdot \left( B_{n}U_{n}\right)
-\nabla \cdot \left[ (B_{n}-B)u\right] +F_{n} & \mathrm{in}\ Q, \\
{\partial }_{\nu }U_{n}=0 & \mathrm{on}\ \Sigma , \\
U_{n}(x,0)=0\  & x\in \Omega .%
\end{cases}
\label{ynzncap}
\end{equation}%
Multiply the first equation of (\ref{ynzncap}) by $U_{n}$, and integrate
over $\Omega $ to give
\begin{equation*}
\frac{d}{dt}\left\vert U_{n}\right\vert _{2}^{2}+\left\vert \nabla
U_{n}\right\vert _{2}^{2}\leq C\left\Vert B_{n}\right\Vert _{\infty
}^{2}\left\vert U_{n}\right\vert _{2}^{2}+C\left\Vert u\right\Vert _{\infty
}^{2}\left\vert B_{n}-B\right\vert _{2}^{2}+C\int_{\Omega }F_{n}Y_{n}dx.
\end{equation*}%
It then follows from Gronwall's lemma that
\begin{equation}
\left\vert U_{n}(\cdot ,T)\right\vert _{2}^{2}\leq e^{C\left\Vert
B_{n}\right\Vert _{\infty }^{2}T}\left( C\left\Vert u\right\Vert _{\infty
}^{2}\left\Vert B_{n}-B\right\Vert _{2}^{2}+C\int_{\Omega
}F_{n}Y_{n}dx\right) .  \label{ynzncapp}
\end{equation}%
By \eqref{fcn1}, \eqref{ynbn}, and \eqref{ynbn1}, we see that the right-hand
side of \eqref{ynzncapp} tends to $0$ as $n\rightarrow \infty $, and hence $%
\left\vert U_{n}(\cdot ,T)\right\vert _{2}\rightarrow 0$. Since $%
u_{n}(x,T)=0 $ for almost all $x\in \Omega $, we get that $u(x,T)=0$ for
almost all $x\in \Omega $. It then follows that $u\in \Psi (B)=\Lambda (\xi
) $. Therefore, $\Lambda $ is upper semi-continuous.

Finally, it remains to show that $\Lambda (K)\subset K$. Indeed, by the
standard energy estimate we see that for any $\xi \in K,$ each element $u$
of $\Lambda (\xi )$ satisfies
\begin{equation*}
\left\Vert u\right\Vert _{L^{\infty }(0,T;L^{p}(\Omega ))}\leq
e^{C(1+\left\Vert B\right\Vert _{\infty }^{2})(1+T)}\left( \left\vert
u_{0}\right\vert _{p}+\left\Vert \mathbf{1}_{\omega }f\right\Vert
_{p}\right) ,
\end{equation*}%
which together with \eqref{kkk0} in Proposition \ref{lemmap} and %
\eqref{fcontrolin0} leads to $u\in L^{\infty }(Q)\cap L^{\infty
}(0,T;L^{p}(\Omega ))\ $and
\begin{equation*}
\left\Vert u\right\Vert _{L^{\infty }(0,T;L^{p}(\Omega ))}+\left\Vert
u\right\Vert _{\infty }\leq e^{c_{1}\kappa _{0}}\left\vert u_{0}\right\vert
_{\infty },
\end{equation*}%
where $c_{1}$ is a positive constant independent of $T$, and $\kappa_0$ is given by \eqref{k0}. If $\left\vert u_{0}\right\vert
_{\infty }\leq e^{-c_{1}\kappa _{0}}$ which is exactly \eqref{intialdata},
then $\left\Vert u\right\Vert _{\infty }\leq 1$ and hence $\Lambda
(K)\subset K$. Apply Kakutani's fixed point theorem to obtain at least one
fixed point $u$ of $\Lambda $; that is $u\in \Lambda (u)$. This $u$ together
with $v=v_{u}$, the solution of \eqref{lin-e} with $\xi =u$, gives the
solution of \eqref{e}, corresponding with some control $f$ and $u(x,T)\equiv
0$. This completes the proof.\hfill $\Box $

\end{document}